\documentclass[12pt]{article}

\title{Logarithms Over a Real Associative Algebra}

\author{Nathan BeDell \\ nbedell@tulane.edu}

\usepackage[a4paper, total={6in, 8in}]{geometry}
\usepackage{graphicx}
\usepackage{amsmath}
\usepackage{amssymb}
\usepackage{mathrsfs}
\usepackage{amsfonts}
\usepackage[utf8]{inputenc}
\usepackage[english]{babel}
\usepackage{amsthm}
\usepackage{amscd}
\usepackage{tikz}
\usepackage{bbm}

\newcommand{\Acal}{\mathcal{A}}

\newcommand{\Ld}{\mathrm{Ld}}

\theoremstyle{definition}
\newtheorem{definition}{Definition}[section]

\newtheorem{theorem}[definition]{Theorem}
\newtheorem{example}[definition]{Example}
\newtheorem{lemma}[definition]{Lemma}
\newtheorem{corollary}[definition]{Corollary}
\newtheorem{proposition}[definition]{Proposition}

\theoremstyle{remark}

\theoremstyle{remark}

\theoremstyle{remark}
\newtheorem*{procedure}{Procedure}

\theoremstyle{remark}

\begin{document}

\maketitle

\begin{abstract}
Extending the work of Freese \cite{Freese15} and Cook \cite{Cook12}, which develop the basic theory of calculus and power series over real associative algebras, we examine what can be said about the logarithmic functions over an algebra. In particular, we find that for any multiplicative unital nil algebra the exponential function is injective, and hence the algebra has a unique logarithm on the image of the exponential. We extend this result to show that for a large class of algebras, the logarithms behave incredibly similarly to the logarithms over the real and complex numbers depending on if they are ``Type-R'' or ``Type-C'' algebras. 
\end{abstract}

\section{Introduction}

The theory of analysis over commutative unital associative algebras, and more generally the theory of analysis over more general contexts such as for non-commutative or even non-associative algebras is long and complex, with authors often re-inventing the wheel, and a number of different approaches which are not all compatible. Surprisingly however, it does not seem that there is any literature available that studies the properties of the logarithm for completely general commutative associative algebras of finite dimension. The closest to this might be \cite{Lorch42}, where the properties of the logarithm are deduced for algebras that take $\mathbb{C}$ as their base field.

In this paper, we wish to develop the theory of the logarithm in the more general context where we take algebras $\mathcal{A}$ with $\mathbb{R}$ as a base field, since an algebra over the complex numbers as a base field can simply be viewed as the real algebra $\mathbb{C} \otimes \mathcal{A}$, since $\mathbb{C}$ is itself a two dimensional $\mathbb{R}$-algebra. 

Throughout this paper, unless we explicitly mention it, by an \textit{algebra} we mean an associative finite dimensional commutative unital algebra over the reals.

\begin{definition}
An algebra $\mathcal{A}$ is a finite dimensional real vector space together with a bilinear multiplication operation $\star : \mathcal{A} \times \mathcal{A} \rightarrow \mathcal{A}$ satisfying the following properties:
\begin{enumerate}
\item $ v \star (w \star z) = (v \star w) \star z$ for all $v,w,z \in \mathcal{A}$
\item There exists an element $\mathbbm{1} \in \mathcal{A}$ such that $\mathbbm{1} \star z = z \star \mathbbm{1} = z$ for all $z \in \mathcal{A}$.
\end{enumerate}
\end{definition}

\noindent If $v \star w = w \star v$ for all $v,w \in \Acal$ then we say $\Acal$ is a commutative algebra. A linear map $\phi : \mathcal{A} \rightarrow \mathcal{B}$ between two algebras is an algebra morphism if it satisfies the homomorphism property $\phi(z \star w) = \phi(z)\star\phi(w)$ and $\phi(\mathbbm{1})=\mathbbm{1}$. 

\begin{proposition}
\label{prop:homo_basis_elements}
Given an algebra $\mathcal{A}$ with basis $\beta = \{ v_1, \dots , v_n \}$ and a linear map $\phi : \mathcal{A} \rightarrow \mathcal{B}$ between two algebras, if $\phi(v_i \star v_j) = \phi(v_i) \star \phi(v_j)$ for all basis elements $v_i, v_j$ and $\phi(\mathbbm{1}) = \mathbbm{1}$, then $\phi(v \star w) = \phi(v) \star \phi(w)$ for all $v,w \in \mathcal{A}$.
\end{proposition}

In addition to this, it will be important to note throughout the paper the following as a consequence of the fact that our algebras are finite dimensional vectors spaces:

\begin{proposition}
\label{prop:algebra_noetherian}
Let $\mathcal{A}$ be a commutative algebra, then $\mathcal{A}$ is a Noetherian and Artinian ring, In particular, this implies that every ideal $I$ of $\mathcal{A}$ is finitely generated.
\end{proposition}

One common way of describing an algebra is by representing it as a suitable subspace of $\mathbb{R}^{n \times n}$, with the usual matrix addition and multiplication representing the multiplication and addition in the algebra. Given a basis $\beta = \{v_1, \dots, v_n\}$ for the algebra, we define a \textit{matrix regular representation}:

\begin{definition}
Given an algebra $\mathcal{A}$ the set of all linear transformations $T: \Acal \rightarrow \Acal$ for which $T(x \star y) = T(x) \star y$ forms the \textit{regular representation} which we denote $\mathcal{R}_{\Acal}$. Clearly $T(x) = T(1 \star x) = T(1) \star x$ hence the regular representation is formed by left-multiplications of $\Acal$.  \\

Denote the left-multiplication by $\alpha \in \Acal$ by $L_\alpha(x) = \alpha \star x$ for each $x \in \Acal$. Given basis $\beta = \{ v_1, \dots, v_n \}$, 
we define $M_\beta(\alpha) = [L_\alpha]_\beta$ and denote the collection of all such matrices by $M_\beta(\mathcal{A})$. 
\end{definition}

Since $\Acal$ is unital and finite dimensional it is well-known that $\Acal$, $\mathcal{R}_{\Acal}$ and $M_\beta(\mathcal{A})$ are isomorphic as algebras. In particular, we have the identifications:
$$ \alpha \ \leftrightarrow \ L_{\alpha} \ \leftrightarrow \ [L_{\alpha}]_{\beta} $$
given a choice of basis $\beta$. For further discussion and some elementary proofs see \cite{Cook12}. Below we illustrate how to calculate $M_\beta(z)$ for an arbitrary element $z \in \mathcal{A}$.

\begin{example}\label{ex:1}
With respect to the basis $\beta = \{1,j\}$, let $z = x + y j$ be an arbitrary element of $\mathcal{H}$. $(x + y j) \star 1 = x + y j$ so the first column of $M_\beta(z)$ is $(x,y)$. Also, $(x + y j) \star j= y + x j$, so the second column of $M_\beta(z)$ is $(y,x)$. Hence: $$ M_\beta(z) = \begin{bmatrix} x & y \\ y & x \end{bmatrix} $$
\end{example}

\begin{example}\label{ex:2}
Let $\beta = \{ 1, i \}$ and let $x = x + i y$ be an arbitrary element of $\mathbb{C}$. Applying the same method as we did in the preceding example, we find $(x + i y) \star 1 = x + i y$ and $(x + i y) \star i = -y + i x$. Hence: $$ M_\beta(z) = \begin{bmatrix} x & -y \\ y & x \end{bmatrix} $$
\end{example}

 Examples \ref{ex:1} and \ref{ex:2} illustrate an interesting result on the structure of matrix regular representations; given a unital basis $\{1, v_2, \dots, v_n\}$ for the algebra $\mathcal{A}$, the $i$th column of the regular representation is an $\mathcal{A}$-multiple of the first column: 

Although this paper is not focused on the calculus of logarithms over an algebra, to keep the paper self-contained we do need to briefly discuss the notion of $\mathcal{A}$-differentiability, and some of the basic properties of the $\mathcal{A}$-derivative that we will use in Section \ref{sec:logarithms}.  For a more complete exposition of this we direct the reader to \cite{Cook12}.

\begin{definition}
Given a function $f : \mathcal{A} \rightarrow \mathcal{A}$, we say that $f$ is $\mathcal{A}$-differentiable at $z_0$ if the Frechet differential $df_{z_0}$ exists, and $df_{z_o} \in \mathcal{R}_{\Acal}$. 
\end{definition}

 This condition implies that $df_{z_0}(z) = L_\alpha(z)$ for some $\alpha \in \mathcal{A}$, and hence allows us to define the \textit{derivative} by setting $f'(z_0) = \alpha$, where $\alpha$ of course in general depends on $z_0$. We also will use the notation $\frac{df(z)}{dz} = f'(z)$ for the $\mathcal{A}$-derivative.

This notion of differentiability satisfies all of the basic properties of the usual real derivative, for example:

\begin{proposition}
\label{prop:derivative_properties}  For $\mathcal{A}$-differentiable functions $f,g : \mathcal{A} \rightarrow \mathcal{A}$,
\begin{enumerate}
\item $(f(g(z)))' = f'(g(z)) \star g'(z)$.
\item $\frac{d}{dz} z^n = n z^n$ for powers $n \in \mathbb{N}$.
\item $(f(z)\star g(z))' = f'(z) \star g(z) + f(z) \star g'(z)$ (here we assume $\Acal$ is commutative)
\item if $c \in \mathcal{A}$ a constant then $(c f(z) + g(z))' = c f'(z) + g'(z)$.
\end{enumerate}
\end{proposition}

\subsection{Semisimple and Nil Algebras}
\label{sec:Algebra_Presentations}

Another convenient way to describe a finite dimensional, commutative, and unital algebra is as the quotient of some real polynomial ring by an ideal. If we have such an algebra $\mathcal{A}$ isomorphic to $\mathcal{P} = \mathbb{R}[x_1, \dots , x_k]/I$ for some $k \in \mathbb{N}$ and $I$ an ideal of $\mathbb{R}[x_1, \dots, x_k]$ we say that $\mathcal{P}$ is a \textit{presentation} of the algebra $\mathcal{A}$. 

\begin{definition}[Standard presentations of typical algebras] $ $
\label{def:common_algebras}
\begin{enumerate}
\item The $n$-hyperbolic numbers: $\mathcal{H}_n := \mathbb{R}[j]/\langle j^n - 1 \rangle$
\item The $n$-complicated numbers: $\mathcal{C}_n := \mathbb{R}[i]/\langle i^n + 1 \rangle$
\item The $n$-nil numbers: $ \mathbf{\Gamma}_n := \mathbb{R}[\epsilon]/\langle \epsilon^n \rangle$
\item The total $n$-nil numbers: $\mathbf{\Xi}_n := \mathbb{R}[\epsilon_1,\dots,\epsilon_n]/\langle \epsilon_i\epsilon_j | i,j \in \{ 1, 2, \dots , n \} \rangle$
\end{enumerate}
\end{definition}

\noindent For example, $\mathcal{C}_2$ is just the usual complex numbers, denoted simply $\mathbb{C}$, and $\mathcal{H}_2$ is just the hyperbolic numbers, denoted simply $\mathcal{H}$. Similarly, we take the convention that $\mathbf{\Gamma}$ by itself denotes $\mathbf{\Gamma}_2$. 

The nil numbers are a special class of what in our terminology we will call \textit{unital nil algebras} -- that is, an algebra with basis $\{1, \epsilon_1, \dots, \epsilon_{n-1} \}$, where each $\epsilon_k$ is \textit{nilpotent}. In other words, for each $\epsilon_k$ there exists $m \in \mathbb{N}$ such that $(\epsilon_k)^m = 0$. This terminology is inspired by the use of  the term \textit{nil algebra} used by Abian \cite{abian71} to refer to an algebra in which every element of the algebra is nilpotent. However, we will sometimes say simply ``nil algebra'' in this paper when we mean ``unital nil algebra'', since in our context all algebras are unital. A unital nil algebra is in some sense the closest you can get to a nil algebra while still being a unital algebra. 

\begin{definition}
\label{def:nil_star}
Given an algebra $\mathcal{A}$, let $\mathrm{Nil}^*(\mathcal{A})$ denote the smallest unital sub-algebra of $\mathcal{A}$ that contains the nilradical of $\mathcal{A}$, $\mathrm{Nil}(\mathcal{A})$ -- that is, the ideal formed from all nilpotent elements of $\mathcal{A}$.
\end{definition}

\begin{proposition}\label{prop:nil_local}
An algebra $\mathcal{A}$ is a unital nil algebra if and only if it is the smallest unital subalgebra of $\mathcal{A}$ containing the nilradical of $\mathcal{A}$. In other words, $\mathcal{A}$ is a unital nil algebra if and only if $\mathcal{A} = \mathrm{Nil}^*(\mathcal{A})$.

\begin{proof}
If $\mathcal{A}$ is a unital nil algebra with unital nil basis $\{1, \epsilon_1, \dots, \epsilon_n \}$, then clearly $\mathrm{Nil}(\mathcal{A}) = \langle \epsilon_1, \dots, \epsilon_n \rangle$, and $\mathcal{A}/\mathrm{Nil}(\mathcal{A}) \cong \mathbb{R}$, so $\mathrm{Nil}(\mathcal{A})$ is maximal, and hence $\mathcal{A}$ is the smallest unital nil algebra containing $\mathrm{Nil}(\mathcal{A})$.

Conversely, suppose that $\mathcal{A} = \mathrm{Nil}^*(\mathcal{A})$, and let $\{ \epsilon_1, \epsilon_2, \dots, \epsilon_n \}$ be a basis for $\mathrm{Nil}(\mathcal{A})$, then clearly $\{ 1, \epsilon_1, \dots, \epsilon_n \}$ is a linearly independent set that spans a subset of $\mathrm{Nil}^*(\mathcal{A})$. Also, we argue that this set must span $\mathrm{Nil}^*(\mathcal{A})$, since it contains all of $\mathrm{Nil}(\mathcal{A})$, and the unit $1$, so by the minimality condition, this must be all of $\mathrm{Nil}^*(\mathcal{A})$. Thus, $\mathcal{A} = \mathrm{Nil}^*(\mathcal{A})$ is a multiplicative nil algebra.
\end{proof}
\end{proposition}

In addition to these basic families, we should also mention the so called $n$-complex numbers $\mathbb{C}_n = \mathbb{C}^{\otimes n}$, where $X^{\otimes n}$ denotes the $n$-fold tensor product of rings.\footnote{For the unfamiliar reader, the tensor product of algebras can be thought of as the algebra which combines the set of generators and relations for a presentation for the algebra. In other words, $\mathbb{R}[x_1, \dots, x_n]/I \otimes \mathbb{R}[y_1, \dots, y_m]/J \cong \mathbb{R}[x_1,\dots,x_n,y_1,\dots,y_m]/(I+J)$. So $\mathbb{C}_2 = \mathbb{C} \otimes \mathbb{C} \cong \mathbb{R}[i_1,i_2]/\langle i_1^2 + 1, i_2^2 +1 \rangle$ for example.} In particular, the analysis of the bicomplex numbers $\mathbb{C}_2$ has been studied extensively, for example, by Price \cite{Price}.

If $\Acal$ is a real vector space with basis $\beta = \{ v_1, \dots , v_n \}$ then given appropriate \textit{structure constants} $c_{ij}^k \in \mathbb{R}$ we may define a multiplication on $\Acal$. In particular, define
$$ v_i \star v_j = \sum_{k=1}^n c_{ij}^k v_k $$ on basis elements, and extend bilinearly to define $\star$ on $\Acal$. Naturally, the structure constants must be given such that the defined multiplication is associative and unital. That said, we typically begin with a given algebra $\Acal$ and simply use the structure constants with respect to a given basis to study the structure of $\Acal$. This notion of structure constants will be important in our derivation of the injectivity of the exponential function in unital nil algebras.

In some sense, the introduction of nilpotent elements into an algebra introduces complications into our theory. Essentially, this is because otherwise, the classification of algebras becomes rather simple.

\begin{definition}
	An algebra is called semisimple if it has no nilpotent elements.
\end{definition}

Usually, one defines a ring to be semisimple if it's Jacobson radical is zero. However, because we are working in the finite dimensional context, our characterization is equivalent to the one involving the Jacobson radical.

\begin{theorem}
	Every commutative semisimple algebra $\Acal$ is isomorphic to $\mathbb{R}^n \times \mathbb{C}^m$ for some $n,m \in \mathbb{N}$.
\end{theorem}

Thus, by a more optimistic characterization, the non-semisimple case, where an algebra has non-trivial nilpotent elements, is where our theory truly departs from the theory of real and complex analysis non-trivially. 

\section{The Exponential and Logarithm Functions for an Algebra}
\label{sec:logarithms}

The exponential function over a commutative algebra always exists, and is defined via power series in $\Acal$:

\begin{definition}
Given an algebra $\mathcal{A}$, we define the exponential function, denoted either $e^z$ or $\exp(z)$ on an algebra by setting $$ e^z = \sum_{k=0}^{\infty} \frac{z^k}{k!} $$
\end{definition}

\noindent Convergence of power series in an algebra is studied in \cite{Freese15} where it is shown:

\begin{theorem}
For any algebra $\mathcal{A}$, the series expansion for the exponential function has radius of convergence $R = \infty$, and hence is well-defined on the whole algebra.
\end{theorem}

 The exponential function for an algebra obeys all of the familiar algebraic properties of the real exponential, in particular:

\begin{proposition}
\label{prop:exponential_properties}
For all $z,w \in \mathcal{A}$:
\begin{enumerate}
\item $e^{0} = 1$
\item $e^{z+w} = e^z e^w$
\item $e^{-z} = \frac{1}{e^z}$
\item $e^z$ is a unit
\end{enumerate}
\end{proposition}

\subsection{Logarithms over Semisimple Algebras}

 We now turn our attention to logarithms over a commutative algebra. (From this point on, an ``algebra" always means a commutative one) The properties of the real and complex logarithms are of course well-studied. Some authors, such as \cite{LorentzNumbers}, have even developed the theory of logarithms over algebras such as the hyperbolic numbers, so a natural question, given the results of Freese's series methods that the exponential function exists in any algebra is what can be said about the logarithm in an arbitrary algebra.

By a logarithm, we mean an inverse function to $e^z$ for an algebra. Thus, for our convenience, we use the notation $\mathrm{Ld}(\mathcal{A}) = \mathrm{Im}(exp)$ to denote the \textit{logarithmic domain}, or equivalently, the image of the exponential for an algebra. Specifically, since the exponential function may not be injective (for example, as is the case with $e^z : \mathbb{C} \rightarrow \mathbb{C}$), given an algebra $\mathcal{A}$, and a connected subset $B \subseteq \mathcal{A}$ of the algebra such that $e^z \lvert_B$ is injective, and $\mathrm{Im}(e^z \lvert_B) = \mathrm{Ld}(\mathcal{A})$ we wish to find a function $\log : \mathrm{Im}(\exp) \rightarrow B$ such that $\log(e^z) = z$ for all $z \in B$, and $e^{\log(z)} = z$ for all $z \in \mathrm{Ld}(\mathcal{A})$. We call this a \textit{branch} of the logarithm for the algebra, borrowing some terminology from complex analysis. 

From this, we can define logarithms to other bases, which we denote by $\log_b(z)$, using the standard change of basis formula for logarithms, namely $\log_b(z) = \frac{\log(z)}{\log(b)}$. In addition to this, and a fact that is more relevant to our unpublished work concerning ODEs over an algebra, we may also from the existence of a logarithm define arbitrary power functions over an algebra, i.e. by defining $a^b = e^{b \log(a)}$ for $a, b \in \mathcal{A}$.

We prove first the existence of a logarithm on at least a sub-domain of the image of the exponential for an arbitrary algebra using series methods:

\begin{theorem}
Every algebra $\mathcal{A}$ has a function $\exp^{-1} 
: V \subseteq \Ld(\Acal) \rightarrow \mathcal{A}$ which is inverse to the exponential function $\exp : U \subseteq \mathcal{A} \rightarrow \mathcal{A}$.

\begin{proof}
Begin by defining 

$$ f(z) = \sum_{k=1}^{\infty} \frac{(-1)^{k-1}}{k} z^k $$ and $$ g(z) = \sum_{k=1}^{\infty} \frac{1}{k!}z^k$$ 

Notice then that $g(z) = e^z - 1$, and to align with the standard series definition of $\log$, we define in our algebra context $log(z) = f(z-1)$, so that $f(z) = log(z+1)$. 

\noindent By termwise differentiation of the series (see \cite{Freese15}), we obtain:

$$ f'(z) = \sum_{k=1}^{\infty} (-1)^{k-1} z^{k-1} = \frac{1}{1+z} $$

\noindent and 

$$ g'(z) = \sum_{j=1}^{\infty} \frac{j}{j!} z^{j-1} = \sum_{k=0}^\infty \frac{z^k}{k!} = 1 + g(z) $$

 and hence, if we define $p(z) = f(g(z))$, by the chain rule we obtain $p'(z) = f'(g(z))g'(z) = \frac{1}{1+g(z)} (1 + g(z)) = 1$, since for any $z \in \mathcal{A}$, $e^z = 1 + g(z) \in U(\mathcal{A})$, and hence $\frac{1}{1+g(z)}$ is well defined in the algebra.

\noindent Finally, since $p(0) = 0$, we can conclude that $p(z) = f(g(z)) = z$ for all $z \in \mathcal{A}$.
\end{proof}

\end{theorem}

Although this existence theorem is useful for us, there is more to say about the specifics of the possible branches of the logarithm over an algebra, but first we must give some new terminology for algebras: 

\begin{definition}
\label{def:typeRC}
Let $\mathcal{A}$ be a commutative algebra. If there exists another algebra $\mathcal{B}$ such that $\mathcal{A} \cong \mathbb{C} \otimes \mathcal{B}$, then we say $\mathcal{A}$ is a type-C algebra. 

If $\mathcal{A}$ is not isomorphic to the direct product of an algebra $\mathcal{B}$ and some type-C algebra $\mathcal{C}$, then we say that $\mathcal{A}$ is a type-R algebra.
\end{definition}

Essentially, a type-C algebra is a algebra which can be viewed as an algebra with $\mathbb{C}$ as its base field. A type-R algebra then, is an ``honest to goodness" algebra over $\mathbb{R}$ in the sense that not even part of the algebra can be viewed as being built over $\mathbb{C}$. Our main results in this section revolve around trying to describe the properties of the logarithm in an algebra based on whether it is a type-R, or a type-C algebra.

Notice however, that given Definition \ref{def:typeRC}, not every algebra is either type-R or type-C. For a basic example, consider $\mathcal{H}_3 \cong \mathbb{R} \times \mathbb{C}$, which has both both a ``real piece" and a ``complex piece", and thus is neither type-R nor type-C. It follows trivially from the definition and by the finite-dimensional nature of our algebras that this is always the case.

In order to prove that the logarithms over an algebra behave as expected over type-R and type-C algebra, we will first recall some preliminary propositions and definitions from ring theory.

\begin{definition}
	A local (commutative) ring is a ring with exactly one maximal ideal. We say that an algebra is a local is it is local as a ring.
\end{definition}

\begin{proposition}
	Any finite dimensional algebra is the finite product of finite dimensional local algebras.
\end{proposition}

\noindent With which we prove the following classification theorem:

\begin{theorem}
	\label{thm:local_ring_classification_thm}
	Any algebra is the product of a finite number of copies of $\mathbb{R}$, $\mathbb{C}$, and a finite number of different untial nil algebras.
	
	\begin{proof}
		By the preceding proposition, every finite dimensional algebra is the finite product of local rings. Also, if a ring is local, its unique maximal ideal is identical to its Jacobson radical, which in our context is the same thing as $\mathrm{Nil}(\Acal)$, and hence, by Proposition \ref{prop:nil_local}, the local finite dimensional algebras are exactly the same as the unital nil algebras\footnote{While we tend to think of unital nil algebras as algebras such as $\Gamma$ with non-trivial nil radicals, $\mathbb{R}$ and $\mathbb{C}$ technically also satisfy our definition of unital nil algebra.}.
	\end{proof}
\end{theorem}

 Now that we have this result, we are ready to begin our investigation of logarithms over semisimple algebras, which is rather simple with the help of Wedderburn's classification theorem. 

\begin{lemma}
\label{lem:direct_prod_exponential}
Let $\mathcal{A} = \mathcal{A}_1 \times \dots \times \mathcal{A}_k$, and $e^z : \mathcal{A} \rightarrow \mathcal{A}$ defined as usual, where $z = e_1 x_1 + \dots + e_k x_k$, and $x_1, \dots , x_k$ are elements of the algebras $\mathcal{A}_1, \dots , \mathcal{A}_k$ respectively, then $e^z = e_1 e^{x_1} + e_2 e^{x_2} + \dots + e_k e^{x_k}$.

\begin{proof}
Applying the definition of termwise multiplication in the direct product algebra $ \mathcal{A}_1 \times \dots \times \mathcal{A}_k$ we obtain:
\begin{equation*}
\begin{aligned}
e^z &= e^{e_1 x_1 + e_2 x_2 + \dots + e_k x_k} \\ &= 1 + (e_1 x_1 + e_2 x_2 + \dots + e_k x_k) + \frac{(e_1 x_1 + e_2 x_2 + \dots + e_k x_k)^2}{2!} + \dots \\ &=  1 + (e_1 x_1 + e_2 x_2 + \dots + e_k x_k) + \frac{(e_1 x_1^2 + e_2 x_2^2 + \dots + e_k x_k^2)}{2!} + \dots
\end{aligned}
\end{equation*}

 And hence, since in the direct product algebra $1 = e_1 + e_2 + \dots + e_k$, collecting the power series component-wise yields: $$ e^z = e_1\left(1 + x_1 + \frac{x_1^2}{2!} + \dots\right) + \dots + e_k\left(1 + x_k + \frac{x_k^2}{2!} + \dots\right) = e_1 e^{x_1} + \dots + e_k e^{x_k} $$ Thus completing the proof.

\end{proof}
\end{lemma}

\begin{theorem}
Let $\mathcal{A} = \mathcal{A}_1 \times \dots \times \mathcal{A}_k$, then if $log_1(x_1), log_2(x_2), \dots , log_k(x_k)$ denote branches of the inverses to the exponential functions in the respective algebras $\mathcal{A}_1, \dots , \mathcal{A}_k$ where $B_1, \dots, B_k$ are the images of said branches for which $log_1(x_1)$, $\dots$, $log_k(x_k)$ are inverses to their respective exponential functions, then the function $log(z) = (log_1,log_2, \dots , log_k) : \mathrm{Ld}(\mathcal{A}_1) \times \mathrm{Ld}(\mathcal{A}_2) \times \dots \times \mathrm{Ld}(\mathcal{A}_k) \rightarrow B_1 \times B_2 \times \dots \times B_k $ is the inverse function of $e^z : \mathcal{A} \rightarrow \mathcal{A}$ on the branch $B_1 \times B_2 \times \dots \times B_k$.

\begin{proof}
By Lemma \ref{lem:direct_prod_exponential}, if $z = e_1 x_1 + \dots e_k x_k$, and $x_1, \dots , x_k$ are elements of the algebras $\mathcal{A}_1, \dots , \mathcal{A}_k$ respectively, then $e^z = e_1 e^{x_1} + e_2 e^{x_2} + \dots + e_k e^{x_k}$, and hence by definition $(\log \circ \exp)(z) = e_1 \log_1(e^{x_1}) + \dots + e_k \log_k(e^{x_k}) = z$ for all $z \in B$, and $(\exp \circ \log)(z) = e_1 e^{log_1(x_1)} + \dots + e_k e^{log_k(x_k)} = z$ for all $z \in \mathrm{Ld}(\mathcal{A}_1) \times \mathrm{Ld}(\mathcal{A}_2) \times \dots \times \mathrm{Ld}(\mathcal{A}_k)$.
\end{proof}
\end{theorem}

\begin{corollary}
\label{cor:log_branches_semisimple}
Given a choice of branch cut for each of the complex components, the logarithm for a semisimple algebra $\mathcal{A}$ may be induced from the product of functions $(\log,\dots ,\log, \mathrm{Log}_1, \dots, \mathrm{Log}_m) : (\mathbb{R}^+)^n \times (\mathbb{C}^\times)^m \rightarrow \mathbb{R}^n \times \mathbb{C}^m$ under the isomorphism $\phi : \mathcal{A} \rightarrow \mathbb{R}^n \times \mathbb{C}^m$ given by Wedderburn's theorem, where $\mathrm{Log}_i : \mathbb{C}^\times \rightarrow \mathbb{C}$ denote the chosen branch cut of the complex logarithm.\footnote{Note that usually a branch cut of the complex logarithm is taken with a domain of $\mathbb{C}^\alpha$ -- that is, $\mathbb{C}$ with a ray starting from 0 removed. This is done to ensure that the logarithm is continuous, but sacrifices having an inverse on all of $\mathrm{Ld}(\mathbb{C}) = \mathbb{C}^\times$, however this corollary can be modified to better suit either approach to branch cuts of the complex logarithm.}
\end{corollary}

 This theorem not only gives us a simple way to construct the logarithm for a commutative semisimple algebra from the basic real and complex logarithms we know and love without using series methods, but more importantly allows us to read off the possible domains of the logarithms for an algebra after passing the branches through the isomorphism map to $\mathbb{R}^n \times \mathbb{C}^m$.

\subsection{Nil Exponentials}
\label{sec:nil_exponentials}

Our goal in this section is to show that the exponential function in a unital nil algebra is always injective. We will then prove that logarithms behave as expected over type-R and type-C algebras, using Lemma \ref{lem:direct_prod_exponential} to extend our result for semisimple algebras to arbitrary algebras once we understand how to deal with the nilpotent pieces.

\noindent Consider the following ``generalized divisibility'' ordering :

\begin{definition}
	Let $\beta = \{ 1, v\_2 , \dots, v_n \}$ be a unital basis for $\Acal$. Set $v_i \preceq v_j$ if and only if there there exists a $v_k \in \beta$ such that $c_{ij}^k \neq 0$, where $c_{ij}^k$ denotes the structure constants for $\Acal$ with respect to the basis $\beta$.
\end{definition}

\begin{proposition}
	$\preceq$ is a partial order.
	
	\begin{proof}
		To show reflexivity, since $\beta$ is a unital basis, for all $v_i \in \beta$ we simply take $1 = v_1 \in \beta$ so that $v \star v_1 = v$. In other words, in terms of structure coefficients, $c_{1i}^i = 1$ is non-zero, hence $v_i \prec v_i$ for all $i = 1, \dots, n$. 
		
		Now, suppose that $v_i \prec v_j$ and $v_j \prec v_k$. By definition, there exists $v_\alpha \in \beta$ such that $c_{i\alpha}^j \neq 0$ and $v_\gamma \in \beta$ such that $c_{j\gamma}^k \neq 0$. Also, by definition of the structure constants we have $$ v_i \star (v_\alpha \star v_\gamma) = (\sum_{k = 1}^{v} c_{i\alpha}^k v_k) \star v_\gamma = \sum_{k = 1}^n c_{i\alpha}^k v_k v_\gamma  = \sum_{i=1}^{n} c_{i\alpha}^k (\sum_{\ell = 1}^{n} c_{k\gamma}^\ell v_\ell) = \sum_{i=1}^{n} \sum_{\ell = 1}^{n} c_{i\alpha}^k c_{k\gamma}^\ell  v_\ell$$ Hence, either $c_{i\rho}^k$ is non-zero for some $\rho$ or the $v_i \star (v_\alpha \star v_\gamma)$ must be zero, but this cannot be the case by the assumption that $c_{i\alpha}^j$ and $c_{j\gamma}^k$ are non-zero. Thus, $\preceq$ is transitive.
	\end{proof}
\end{proposition}

 Using this ordering, we will show that products of power series of indeterminates times one of the nilpotent basis elements $v_2, \dots, v_n$ are of a very particular form.

\begin{lemma}
	\label{lem:inductive}
	Given an algebra $\Acal$ with basis $\beta = \{1, v_2, \dots, v_n\}$, we will consider the indeterminates $z_1, z_2, \dots, z_n$, and say that $z_i$ is the indeterminate associated with the basis element $v_i$ for each $i$, and let $p$ be a polynomial in $A$ with indeterminates $z_1, \dots, z_n$. Let $P(p)$ be the proposition that for all $i = 1, 2, \dots, n$, the coefficient of $v_i$ is a polynomial in indeterminates $z_j$ such that $v_j \preceq v_i$. Then we have the following:
	
	\begin{enumerate}
		\item $P$ holds for all polynomials of the form $z_iv_i$, which we will call ``atomic'' polynomials.
		\item $P(p)$ and $P(q)$ implies $P(p + q)$.
		\item If $c \in R$, then $P(p)$ implies $P(cp)$.
		\item $P(p)$ and $P(q)$ implies $p(p\star q)$.
	\end{enumerate}

	\begin{proof}
		Part one of the lemma holds trivially, so we first consider part two. Suppose that $F = f_1v_1 + f_2 v_2 + \dots + f_n v_n = F$ and $G = g_1v_1 + g_2 v_2 + \dots + g_n v_n$ both satisfy $P$, where the $f_i, g_j$ are polynomials satisfying the conditions in the lemma, then $$F + G = (f_1 + g_1)v_1 + (f_2 + g_2)v_2 + \dots + (g_n + f_n)v_n$$ clearly also satisfies $P$. Similarly, $cF = (cf_1)v_1 + \dots (cf_n)v_n$ also satisfies $P$. Finally, consider:
		
		\begin{equation*}
		\begin{aligned}
			FG &= \sum_{i,j=1}^n (f_iv_i)(g_jv_j) = \sum_{i,j=1}^n f_ig_j(v_iv_j) \\
			   &= \sum_{i,j=1}^n f_i g_j (\sum_{k=1}^n c_{ij}^k v_k) = \sum_{i,j,k=1}^n (c_{ij}^k f_i g_j) v_k 
		\end{aligned}
		\end{equation*}
		
		 Then, since by definition $c_{ij}^k$ is non-zero if and only if $v_i \preceq v_k$ and $v_j \preceq v_k$, $FG$ satisfies $P$.
	\end{proof}
\end{lemma}

\begin{corollary}
	\label{cor:nil_power_series}
	Let $\Acal$ be a unital nil algebra with nil basis $\beta$, and let $v, w \in \beta$ and $v,w \neq 1$. Furthermore, suppose that $x,y$ are two indeterminates associated with $v$ and $w$ respectively, then if $p,q$ are power series in $xv$ and $yw$ respectively, then $p\star q$ is a polynomial satisfying $P$ from Lemma \ref{lem:inductive}.
	
	\begin{proof}
		If $v, w \neq 1$, then by the definition of a unital nil algebra they must be nilpotent. Thus, the series $p,q$ both truncate to polynomials. Furthermore, since these polynomials are the linear combination of atomic polynomials, by Lemma \ref{lem:inductive} $P$ holds for both $p$ and $q$, and hence the product $p \star q$ again by Lemma \ref{lem:inductive} . 
	\end{proof}
\end{corollary}	

 Unfortunately, there does not seem to be a simple proof that the ordering $\preceq$ is always antisymmetric (and thus a \textit{poset}), since this is necessary for our technique in the following theorem. Thus, until we find a method to prove that $\preceq$ is in fact always antisymmetric, we will restrict our results to the class of \textit{multiplicative nil algebras}.

\begin{definition}
	Let $\mathcal{A}$ be an algebra. We say a basis $\beta = \{ v_1, \dots , v_n \}$ of $\mathcal{A}$ is multiplicative if for all $v_i,v_j \in \beta$ we have $v_i \star v_j = c v_k$ for some $c \in \mathbb{R}$ and $v_k \in \beta$. If an algebra $\mathcal{A}$ admits a multiplicative basis, then we say $\mathcal{A}$ is a multiplicative algebra.
\end{definition}

 It is easy to show from the definition and the properties of unital nil bases that for such algebras, $\prec$ is in fact a poset\footnote{This is what I call, in \cite{BeDell_AlgebraPaper}, the nil poset for an algebra with respect to a basis.}. Finally, in this context we are able to prove that the exponential function in a unital nil algebra is injective:

\begin{theorem}
\label{thm:nil_exponental_injective}
In any multiplicative commutative unital nil algebra $\mathcal{A}$, the exponential function is injective.

\begin{proof}
After choosing a basis, we identify $\exp : \mathcal{A} \rightarrow \mathcal{A}$ as a function $\exp : \mathbb{R}^n \rightarrow \mathbb{R}^n$, then let $(y_1, y_2, \dots , y_n) = f(x_1, \dots x_n)$ be an element of the image of $\exp$. We will show that $\exp$ is injective by giving an algorithm to construct the unique inverse map $\exp^{-1} : \Ld(\Acal) \rightarrow \mathbb{R}^n$. We will accomplish this by first constructing functions $f_1, \dots f_n$ such that $f_i(f(x_1,x_2, \dots , x_n)) = x_i$, then piecing these together, defining $\exp^{-1} = (f_1,\dots,f_n)$ so that: $$ \exp^{-1}(\exp(x_1,x_2,\dots,x_n)) = (f_1(x_1, \dots, x_n), \dots , f_n(x_1, \dots, x_n)) = (x_1,x_2, \dots , x_n) $$ and hence that $\exp^{-1}$ is in fact the inverse function to $\exp$.

Recall from Proposition \ref{prop:exponential_properties} that the standard rules of exponents still hold in an algebra, and hence $e^{x_1 + x_2 v_2 + \dots + x_n v_n} = e^{x_1}e^{x_2 v_2} \dots e^{x_n v_n}$. Also note that the coefficient of the basis element $1$ is $1$ in the product $e^{x_2 v_2} \dots e^{x_n v_n}$ (we will denote this by $\mathrm{coeff}(1) = 1$ from now on for brevity), since each factor has the form $e^{x_i v_i} = 1 + x_i v_i + \frac{1}{2} (x_i v_i)^2 + \dots$, and thus the only possible product yielding some constant multiple of $1$ is in fact simply $1\cdot1\cdot1 \dots \cdot 1 = 1$. Thus, $\mathrm{coeff}(1) = e^{x_1}$ in the product $e^{x_1}e^{x_2 v_2} \dots e^{x_n v_n}$, and hence in our identification with $\mathcal{A}$ as $\mathbb{R}^n$, we may define $f_1(y_1,\dots,y_n) = \log(y_1)$, which gives us $f_1(\exp(x_1, \dots , x_n)) = x_1$.

We will now iterate a procedure on the set of minimal elements of the lattice $\mathcal{N}_\mathcal{A} - E$, in which $E$ is a subset of $\mathcal{N}_\mathcal{A}$ which we initially set $E = \{ 1 \}$, in order to define the remaining component functions $f_2, \dots, f_n$ of $\exp^{-1}$.

For the initial step of this procedure, select a minimal element $v_i$ of $\mathcal{N}_\mathcal{A} - E = \mathcal{N}_\mathcal{A} - \{ 1 \}$. Since the only element below $v_i$ in $\mathcal{N}_\mathcal{A}$ is 1\footnote{This is where the antisymmetry condition for $\preceq$ is important, as otherwise this may not be true}, by Corollary \ref{cor:nil_power_series}, the coefficient of $v_i$ in the product $e^{x_2 v_2} \dots e^{x_n v_n}$ is a polynomial in the indeterminate $x_i$. Also, notice that the product $e^{x_2 v_2} \dots e^{x_n v_n}$ when expanded contains the term $x_iv_i$, and hence the coefficient of $v_i$ in this product has the form $x_i + p(x_i)$ for some polynomial $p$. Furthermore, $p(x_i)$ must be zero, since the only possible product among powers of basis elements $v_1, v_2, \dots , v_n$ yielding a term with coefficient $v_i$ will be the given product of terms containing $x_iv_i$ and $1$. Hence, we set $f_i(y_1,y_2, \dots , y_n) := \frac{y_2}{y_1}$.

This initial procedure may then be iterated on all other $v_i$ minimal in $\mathcal{N}_\mathcal{A}$, after which we set $E := E + \{ v_{i_1}, \dots , v_{i_k} \}$, where $\{ v_{i_1}, v_{i_2}, \dots, v_{i_k} \}$ is the set of elements minimal in $\mathcal{N}_\mathcal{A} - \{1\}$.

The preceding step was the ``base case'' of our algorithm. We now consider the general step of the algorithm, with $E$ set to $\{ 1, v_{i_1}, \dots , v_{i_k} \}$ and $f_{i_1}, \dots f_{i_k}$ correctly defined as we discussed in the outset of the proof. As before, we now consider elements $v_i$ minimal in the poset $\mathcal{N}_\mathcal{A} - E$, and again similarly we argue that the coefficient of $v_i$ in the product $e^{x_2 v_2} \dots e^{x_n v_n}$ must be of the form $y_i = x_i + p(x_{j_1}, \dots , x_{j_l})$, where $p$ is some polynomials in the indeterminates $x_{j_1}, \dots , x_{j_l}$ associated with the elements $v_a$ such that $v_a \preceq v_i$. Now, since every element below $v_i$ is in $E$, we have already constructed the functions $f_{j_1}, f_{j_2}, \dots , f_{j_l}$, so we set
$$f_i(y_1, \dots , y_n) := \frac{1}{y_1} y_i - p(f_{j_1}(y_1, \dots , y_n),f_{j_2}(y_1, \dots, y_n),\dots,f_{j_l}(y_1, \dots, y_n))$$ so that 
\begin{equation*}
\begin{aligned}
 f_i(e^{x_1 + x_2 v_2 + \dots + x_n v_n}) &= \frac{e^{x_1}y_i}{e^{x_1}} - p(f_{j_1}(e^{x_1 + x_2 v_2 + \dots + x_n v_n}), \dots , f_{j_l}(e^{x_1 + x_2 v_2 + \dots + x_n v_n})) \\ &= y_i - p(x_{j_1},x_{j_2}, \dots, x_{j_l}) \\ &= x_i + p(x_{j_1},x_{j_2}, \dots, x_{j_l}) - p(x_{j_1},x_{j_2}, \dots, x_{j_l}) \\ &= x_i
\end{aligned}
\end{equation*}

\noindent and finally, to complete the general iterative step, we set $E := E \cup \{v_i\}$.

\noindent Given our setup, we may now summarize the complete algorithm as follows:

\begin{procedure}
Begin by setting $E = \{ 1 \}$, and iterate the procedure as described in paragraph 4 of the proof, constructing the functions $f_{i_1}, \dots , f_{i_{k_1}}$ and after completing the procedure, setting $E := E \cup \{v_{i_1}, \dots , v_{i{k_1}}\} $.
Next, iterate the procedure described starting in paragraph 6 of the proof, defining $f_{i_{k_1 + 1}}, \dots , f_{i_{k_2}}$, and setting $E := E \cup \{ \}$ at the end of each iteration until $\mathcal{C}_\mathcal{A} - E = \emptyset$, at which point the algorithm terminates, yielding the desired functions $f_1, \dots , f_n$.
\end{procedure}

Since each step of the algorithm described above decreases the size of the finite set $\mathcal{N}_\mathcal{A} - E$, the algorithm is clearly productive, and the algorithm halts whenever $\mathcal{N}_\mathcal{A} - E = \emptyset$, this algorithm terminates giving us $f_1, f_2, \dots , f_n$ such that $\exp^{-1} = (f_1,f_2,\dots,f_n)$ is the inverse of $\exp$.
\end{proof}
\end{theorem}

\begin{example}
In $\mathbf{\Gamma}_3$: 
\begin{align} \notag
e^z &= e^{x+\epsilon y + \epsilon^2 z} \\ \notag
&= e^x e^{\epsilon y} e^{\epsilon^2 z} \\ \notag
&= e^x\left(1 + \epsilon y + \frac{\epsilon^2 y^2}{2}\right)\left(1 + \epsilon^2 z \right) = e^x\left(1 + \epsilon y + \epsilon^2\left(z + \frac{y^2}{2}\right)\right)
\end{align}
And hence, following the algorithm given in Theorem \ref{thm:nil_exponental_injective}, or simply by inspection, we find: $$ \mathrm{Log}_{\mathbf{\Gamma}_3}(x + \epsilon y + \epsilon^2 z) = \log(x) + \epsilon \frac{y}{x} + \epsilon^2 \left(\frac{z}{x} - \frac{y^2}{2x^2}\right) $$
\end{example}

 In addition to this result, it is easy to see from the structure of the component functions of the exponential as deduced in the proof of Lemma \ref{lem:direct_prod_exponential} that:

\begin{proposition}
Given a commutative multiplicative unital nil algebra $\mathcal{A}$ with point set identified as $\mathbb{R}^n$, $\mathrm{Ld}(\mathcal{A}) = \mathbb{R}^+ \times \mathbb{R}^{n-1}$. Similarly, if we consider the complexification $\mathbb{C} \otimes \mathcal{A}$ as identified with $\mathbb{C}^n$, then $\mathrm{Ld}(\mathbb{C} \otimes \mathcal{A}) = \mathbb{C}^\times \times \mathbb{C}^{n-1}$.
\end{proposition}

 Since the component functions of the inverse we constructed to the exponential are only undefined where $\log(x_1)$ is undefined, which is on $\mathbb{R}^- \cup \{ 0 \}$ for the real logarithm, and on $\{ 0 \}$ for the complex logarithm.

\subsection{Final Results}

We now are able to completely characterize the possibilities for logarithms in an arbitrary multiplicative algebra\footnote{It is easy to see that the direct product factors of a multiplicative algebra must be again multiplicative.}. Given the results of Theorem \ref{thm:nil_exponental_injective} we state the following result characterizing the logarithms of any multiplicative algebra, extending the results of Corollary \ref{cor:log_branches_semisimple}:

\begin{proposition}
Given a multiplicative algebra $\mathcal{A}$, then let $\mathcal{A} \cong \mathcal{A}_1 \times \mathcal{A}_2 \times \dots \times \mathcal{A}_k \times \mathbb{R}^n \times \mathbb{C}^m$ be the decomposition of $\mathcal{A}$ proven to exist in Theorem \ref{thm:local_ring_classification_thm}.

If $\mathcal{A}$ is a type-R algebra, then there exists a unique inverse function to $\exp$ on $\Ld(\Acal)$, denoted $\mathrm{Log}_\mathcal{A}(z)$. Otherwise, there exists infinitely many logarithms determined by the branches of the logarithms defined on the complex portion of the algebra.
\end{proposition}

 Thus, at least for the case where our algebra is multiplicative, for a type-R algebra, we see that the exponential function, akin to the real exponential is injective, and hence type-R algebras have unique logarithms defined on all of $\mathrm{Ld}(\mathcal{A})$.

On the other hand, for a type-C algebra, again assuming our conjecture, we have the following result:

\begin{proposition}
\label{prop:type-C-LD}
For a multiplicative type-C algebra $\mathcal{A}$, $\mathrm{Ld}(\mathcal{A}) = \mathcal{A}^\times$
\begin{proof}
Since $(A \times B)^\times = A^\times \times B^\times$, it suffices for us to prove this result on the possible factors of a type C algebra given in Theorem \ref{thm:local_ring_classification_thm}. Also, by complex analysis we already know that $\mathrm{Ld}(\mathcal{A}) = \mathbb{C}^\times$, so it remains to show proposition \ref{prop:type-C-LD} for complexified unital nil algebras. But for such an algebra $\mathcal{A}$ with unital nil basis $\{1,\epsilon_1, \dots, \epsilon_n \}$ so that we may identify the algebra with $\mathbb{C}^n$ we have $\mathrm{Ld}(\mathcal{A}) = \mathbb{C}^\times \times \mathbb{C}^{n-1}$, in other words, in terms of the basis, $\Ld(\Acal) = \{ a + b_1 \epsilon_1 + \dots + b_n \epsilon_n | a,b_1,\dots,b_n \in \mathbb{R}, a \neq 0 \}$ which is precisely the set of units in $\mathcal{A}$.\footnote{This is a consequence of the well-known result from algebra that if $a \in \mathcal{A}^\times$ and $\epsilon, \xi \in \mathrm{Nil}(\mathcal{A})$, then $a + \epsilon \in \mathcal{A}^\times$ and $\epsilon + \xi \in \mathrm{Nil}(\mathcal{A})$}
\end{proof}
\end{proposition}

And hence, although we lose injectivity of the exponential function, we gain the fact that the logarithm is defined on almost all of the algebra. 

\pagebreak
\bibliography{sources}

\end{document}